\documentclass[12pt]{amsart}

\usepackage{amsmath}
\usepackage{amscd}
\usepackage{amsfonts}
\usepackage{amssymb}

\theoremstyle{definition}
\newtheorem{expl}{Example}
\theoremstyle{plain}

\newtheorem{theorem}{Theorem}[section]

\newtheorem{proposition}[theorem]{Proposition}

\newtheorem{corollary}[theorem]{Corollary}

\newcommand{\be}{\begin{equation}}
\newcommand{\ee}{\end{equation}}
\newcommand{\bes}{\begin{equation*}}
\newcommand{\ees}{\end{equation*}}

\newcommand{\cE}{\mathcal{E}}

\newcommand{\cT}{\mathcal{T}}

\begin{document}

\title[Equicontinuity of Operators and Automorphisms]{A Note on Equicontinuity of Families of Operators and Automorphisms}

\author{Orr Shalit}
\date{May 12, 2007}
\begin{abstract}
{This note concerns uniform equicontinuity of families of operators on a separable Hilbert space $H$, and
of families of maps on $B(H)$. It is shown that a one parameter group of automorphisms is uniformly equicontinuous if
and only if the group of unitaries which implements it is so. A ``geometrical" necessary and sufficient condition is given for a family of operators to be uniformly equicontinuous.}
\end{abstract}

\maketitle

\section{Introduction}
In this note, we record a few results about uniform equicontinuity of families of operators on
a Hilbert space $H$ and of families of maps on $B(H)$. The motivation for this study was an attempt to prove
a ``commutant lifting" theorem regarding an E-semigroup and a CP-semigroup that commutes with it (see \cite{SN} for more details). The results that we obtained were actually of a negative type, in the sense that they showed that the methods which we intended to use in order to prove this ``commutant lifting" theorem cannot work. However, we believe that the results obtained in this study are interesting in their own right.

Throughout, $H$ will denote a separable (unless stated otherwise) Hilbert space, $B_1$ will denote the
closed unit ball in $B(H)$, and $H_1$ will denote the closed unit ball in $H$.

Let $\{h_i\}_{i=1}^\infty$ be a dense sequence in $H_1$. The weak topology (and also the $\sigma$-weak topology,
because these two coincide on $B_1$) is induced by the metric
\be\label{eq:d}
d(A,B) = \sum_{i,j=1}^\infty \frac{|\langle(A-B)h_i, h_j \rangle|}{2^{i+j}} \,\, , \,\, A,B \in B_1.
\ee

When we discuss notions of \emph{uniformity} in $B_1$, we shall always mean the uniformity
induced by this metric. A uniformity is also induced on $H_1$ by the metric
\be\label{eq:rho}
\rho(x,y) = \sum_{i=1}^\infty \frac{|\langle x-y, h_i \rangle|}{2^{i}} \,\, , \,\, x,y \in H_1.
\ee
We shall once consider also non-separable $H$, and then the uniformities are given by the connectors
\be\label{eq:connectorH}
\{(x,y) \in H_1 \times H_1 : | \langle x-y,z_i\rangle |<\epsilon\ \, , \, i=1,\ldots, N\}
\ee
for $H_1$, and
\be\label{eq:connectorBH}
\{(A,B) \in B_1 \times B_1 : |\langle (A-B)w_i,z_i\rangle|<\epsilon\ \, , \, i=1,\ldots, N\}
\ee
for $B_1$, where $w_i, z_i \in H_1$, $N\in \mathbb{N}$ and $\epsilon > 0$. When $H$ is separable, these uniformities coincide (that is, (\ref{eq:rho}) with the ones given by the metrics,
because $H_1$ and $B_1$ are compact, and the uniformity for such a space is unique \cite{Wil}.

Recall the notions of \emph{equicontinuity} (henceforth EC) and \emph{uniform equicontinuity} (henceforth UEC).
In the setting of a family $\mathcal{F}$ of maps in a metric space $(X,d)$, the defitions are as follows:
\bes
\textrm{EC  :   } \forall x \in X ,\epsilon >0 . \exists \delta >0 . \forall f\in\mathcal{F} ,y \in X .
d(x,y)<\delta \Rightarrow d(f(x),f(y))<\epsilon \,
\ees
\bes
\textrm{UEC :   } \forall \epsilon >0 . \exists \delta >0 . \forall f\in\mathcal{F} . \forall x,y \in X .
d(x,y)<\delta \Rightarrow d(f(x),f(y))<\epsilon \,
\ees
Let us note that these notions are the same when one considers maps between compact metric spaces, (and probably
also compact Hausdorff uniform spaces (which is what we have if $H$ is not separable)). Here is the proof (for the
metric case):

\begin{proof}
Assume that $\mathcal{F}$ is not UEC. We have sequences $\{f_n\}$ in $\mathcal{F}$ and $\{x_n\},\{y_n\}$ in $X$, and some positive $\epsilon$ such that $d(x_n,y_n)< 1/n$ and $d(f_n(x_n),f_n(y_n))>\epsilon$. We may assume that $x_n\rightarrow x \in X$. But then $y_n \rightarrow x$. Thus
\bes
d(f_n(x),f_n(x_n)) + d(f_n(x),f_n(y_n))\geq d(f_n(x_n),f_n(y_n))>\epsilon ,
\ees
and $\mathcal{F}$ is not EC at $x$.
\end{proof}

\section{The results}
In the following propositions I shall make use of the following easy to verify fact:

\noindent{\bf Fact.}
\emph{If $\mathcal{F}$ and $\mathcal{G}$ are UEC families of a uniform space onto itself, then $\mathcal{F}\circ\mathcal{G}:=\{f\circ g: f\in \mathcal{F},g\in \mathcal{G}\}$ is also UEC.}

\begin{proposition}\label{prop:T}
Let $H$ be an infinite dimensional (not-necessarily separable) Hilbert space. Let $\mathcal{T} \subseteq B_1$. $\mathcal{T}$ is UEC if and only if
the family $\{\psi_T\}_{T\in\mathcal{T}}$ of continuous maps on $B_1$ given by
\bes
\psi_T (A) = TA \,\, , \,\, A\in B_1
\ees
is UEC.
\end{proposition}
\begin{proof}
Assume that $\mathcal{T}$ is UEC. Let $x_i,y_i \in H_1, i=1,\ldots , N,$ and $\epsilon > 0$ be given. To show that $\{\psi_T\}_{T\in\mathcal{T}}$ is UEC, we must find $z_i, w_i \in H_1, i=1, \ldots, M,$ and $\delta > 0$ such that for all $A,B \in B_1$ and $T \in \mathcal{T}$,
\be\label{eq:est}
|\langle(A-B)z_i,w_i\rangle| < \delta \,\, , \,\, i = 1,\ldots,M
\ee
implies
\be\label{eq:est2}
|\langle(TA-TB)x_i,y_i\rangle| < \epsilon \,\, , \,\, i = 1,\ldots,N .
\ee
Let $\delta$ and $u_{j}\in H_1, j=1,\ldots,K$ be such that
$$|\langle \xi-\eta,u_{j}\rangle| < \delta \,\, , \,\, j = 1,\ldots,K $$
implies
$$|\langle T \xi-T \eta,y_i\rangle| < \epsilon \,\, , \,\, i = 1,\ldots,N \, , \, T\in \mathcal{T}.$$
Taking $M = \max \{N,K \}$ and $\{z_i\} = \{x_i\}, \{w_i\} = \{u_j\}$, with repititions if needed, we have that
(\ref{eq:est}) $\Rightarrow$ (\ref{eq:est2}).

Now assume that $\{\psi_T\}_{T\in\mathcal{T}}$ is UEC.
Let $\epsilon>0$ and a unit vector $y\in H_1$ be given. Let $z_i,w_i \in H_1, i=1, \ldots, N$ and $\delta$ be such that
\be\label{eq:imply}
\forall i. |\langle (A-B)z_i,w_i \rangle| < \delta  \Rightarrow  \forall T\in\mathcal{T}.|\langle T(A-B)y,y \rangle| < \epsilon.
\ee
We may assume that, up to repititions, the $z_i$ are linearly independent (by throwing away a few of them and perhaps making $\delta$ smaller). By applying a linear transformation, we may also assume that the $\{z_i\}$ are orthonormal. We distinguish between the two cases $y \in \textrm{span}\{z_1,\ldots,z_N\}$ and $y \notin \textrm{span}\{z_1,\ldots,z_N\}$.

Assume first that $y \in \textrm{span}\{z_1,\ldots,z_N\}$. Since we have already messed with the $z_i$'s, we may as well assume $y = z_1$.
For any $\xi, \eta \in H_1$, we define $A,B \in B_1$ by $A z_1 =  \xi, B z_1 = \eta$ and zero on the complement. Now if
$|\langle \xi - \eta, w_i \rangle| < \delta$ for all $i$, then $|\langle (A-B)z_i,w_i \rangle| < \delta$ for all $i$, from which it follows that for all $T\in\mathcal{T}$,
$$\epsilon > |\langle T(A-B)y,y \rangle| = |\langle T(\xi-\eta),y \rangle| .$$

Assume now that $y \notin \textrm{span}\{z_1,\ldots,z_N\}$. In
this case, we may clearly put $z_{N+1} = y, w_{N+1}=0$, and the
implication in equation (\ref{eq:imply}) still holds (and is not
void), so we are back in the case already considered.
\end{proof}

\begin{proposition}\label{prop:Tstar}
Let $H$ be an infinite dimensional (not-necessarily separable) Hilbert space. Let $\mathcal{T} \subseteq B_1$. $\mathcal{T}^*$ is UEC if and only if
the family $\{\phi_T\}_{T\in\mathcal{T}}$ of continuous maps on $B_1$ given by
\bes
\phi_T (A) = AT \,\, , \,\, A\in B_1
\ees
is UEC.
\end{proposition}
\begin{proof}
The proof follows from the identity
$$\phi_T (A) = (\psi_{T^*}(A^*))^* ,$$
using the fact that $A \mapsto A^*$ is weak operator continuous, hence uniformly weak operator continuous, and the fact stated before the previous proposition.
\end{proof}

\begin{theorem}
Let $\{u_t\}_{t\in\mathbb{R}}$ be a one-parameter unitary group in $B(H)$, and let $\{\alpha_t\}_{t\in\mathbb{R}}$ the group of automorphisms of $B(H)$ implemented by $\{u_t\}_{t\in\mathbb{R}}$. Then $\{u_t\}_{t\in\mathbb{R}}$ is UEC if and only if $\{\alpha_t\}_{t\in\mathbb{R}}$ is UEC.
\end{theorem}
\begin{proof}
If $\{u_t\}_{t\in\mathbb{R}}$ is UEC, then by the above propositions so are the families $\{\phi_{u_t}\}_{t\in\mathbb{R}}$ and $\{\psi_{u_t}\}_{t\in\mathbb{R}}$. Then the family $\{\psi_{u_t} \circ \phi_{u_s}\}_{t,s\in\mathbb{R}}$ is UEC, and in particular its subfamily $\{\alpha_t\}_{t\in\mathbb{R}} = \{\psi_{u_t} \circ \phi_{u_{-t}}\}_{t\in\mathbb{R}}$.

Assume that $\{u_t\}_{t\in\mathbb{R}}$ is not UEC. Then (after making a few selections) there is a sequence $x_n \rightarrow x$ in $H_1$ and a sequence $\{t_n\}$ of real numbers such that $u_{t_n} x_n \rightarrow z$ and $u_{t_n} x \rightarrow w$, with $x\neq z$. We define $A_n = x_n \otimes x_n$ and $A = x \otimes x$. Then
for all $h,g \in H$,
\begin{align*}
\langle A_n h,g \rangle
&= \big\langle \langle h,x_n\rangle x_n,g \big\rangle \\
&= \langle h,x_n\rangle  \langle x_n,g\rangle \\
&\rightarrow \langle h,x\rangle  \langle x,g\rangle \\
&= \big\langle \langle h,x\rangle x,g \big\rangle \\
&= \langle A h,g \rangle
\end{align*}
so $A_n \rightarrow A$. But $\alpha_{t_n} (A_n) = u_{t_n} \circ (x_n \otimes x_n ) \circ u_{t_n}^* = (u_{t_n}x_n) \otimes (u_{t_n}x_n) \rightarrow z\otimes z$, while $\alpha_{t_n} (A) = (u_{t_n}x) \otimes (u_{t_n}x) \rightarrow w \otimes w$. Thus, $\{\alpha_{t}\}_{t\in\mathbb{R}}$ is not EC at $A$.
\end{proof}

So the question of the UEC-ity of a group of automorphisms is equivalent to the question of the UEC-ity of the group of
unitaries implementing it. Is it reasonable to expect that a group of unitaries be UEC? The following example shows that it is not.
\begin{expl}\label{expl:shift}
Let $S$ be the left shift on $\ell^2(\mathbb{Z})$. We reconstruct the metric on the ball of $\ell^2(\mathbb{Z})$ as follows. Let $\{a_n\}$ be some rapidly increasing sequence of integers beginning with $a_1 = 0$. We organize the vectors $h_i$ from equation (\ref{eq:rho}) such that $h_{a_n} = e_n$, and such that
$$\{h_1, h_2, \ldots, h_{a_n}\} \subset \textrm{span}\{e_{-n}, e_{-n+1}, \ldots, e_{n-1}, e_n\} .$$
(See the proof of Theorem \ref{thm:eq})
Now it is clear that $\rho(e_k,e_{k+1})$ can be made arbitrarily small by taking $k$ to be large enough, but
$\rho(S^k e_k,S^k e_{k+1}) \geq 1/2$. This shows that $\{S^k\}_{k\in\mathbb{Z}}$ is not UEC.
\end{expl}
\begin{expl}
Let $u_t$ be given on $L^2[-\pi,\pi]$ by multiplication with the function $e^{itx}$. This is a continuous one-parameter unitary group with a bounded generator (multiplication by $x$). This group is not UEC because the subgroup $\{u_k\}_{k\in\mathbb{Z}}$ is (unitarily equivalent) the shift from Example \ref{expl:shift}.
\end{expl}

In the following theorem we characterize familes of contractions that are UEC on $H_1$.
\begin{theorem}\label{thm:eq}
Let $\mathcal{T} \subset B_1$ be a family of maps on $H_1$. For every finite dimensional subspace $V \subset H$, let $P_V$ be the orthogonal projection on $V$.
The following statements are equivalent:
\begin{enumerate}
    \item For every finite dimensional subspace $V \subset H$, and every $c>0$,
    $$\mathrm{dim}\{x\in V^\perp : \exists T\in\mathcal{T} . \|P_V T x\| \geq c\|x\| \} < \infty .$$
    \item There exists an orthonormal basis $\{e_n\}_{n=1}^\infty$ such that for all $n$, and all $c>0$ $$\mathrm{dim}\{x\in F_n^\perp :  \exists T\in\mathcal{T} . \|P_{F_n} T x\| \geq c\|x\| \} < \infty ,$$
    where $F_n := \mathrm{span}\{e_1,e_2,\ldots,e_n\}$.
    \item $\mathcal{T}$ is UEC on $H_1$.
\end{enumerate}
\end{theorem}
\noindent {\bf Remark.} By $\mathrm{dim}\{x\in V^\perp : \exists T\in\mathcal{T} . \|P_V T x\| \geq c\|x\| \}$ we mean that this set is contained in a finite dimensional
space. It is not a subspace.

\noindent {\bf Remark.} This is quite an indirect way of proving that (2) implies (1).
\begin{proof}
(1) $\Rightarrow$ (2) is true.

(2) $\Rightarrow$ (3). It is sufficient to prove that that $\mathcal{T}$ is EC at every point of $H_1$. Without loss of
generality, we show this for the point $0$. Let $\{h_n\}$ be a (norm) dense sequence in $H_1$ and let $\{a_n\}\subset \mathbb{N}$ be such that $a_n \nearrow \infty$, $e_n = h_{a_n}$ for all $n$, and
$$\mathrm{span}\{h_1, \ldots, h_{a_n}\} \subseteq F_{n} $$
(such sequences can be constructed by constructing $1/n$-nets for $F_n$).
Let $\rho$ be given by equation (\ref{eq:rho}). Let $\epsilon>0$. Take $M$ such that $\sum_{k>a_M} 2^{-k} < \epsilon/4$, and $c<\epsilon/4$. Let $G_0 \subset F_M^\perp$ be a finite dimensional space containing $\{x\in F_M^\perp :  \exists T\in\mathcal{T} . \|P_{F_M} T x\| \geq c\|x\| \}$, and put $G = G_0 \vee F_M$.
Choose $N>M$ such that $\|P_G(I-P_{F_N})\|<\epsilon/4$ (this is possible because $I-P_{F_N}\rightarrow 0$ strongly, and hence in norm on $G$ because $G$ is finite dimensional. Thus $\|(I-P_{F_N})P_G \| \rightarrow 0$, so $\|P_G(I-P_{F_N}) \| \rightarrow 0$ also.). Take $\delta < 2^{-a_N} \cdot \epsilon/4$. We show that for all $x \in H_1$,
$$\rho(x,0) < \delta \Rightarrow \forall T\in\mathcal{T} . \rho(Tx,0) < \epsilon .$$
Assume $x \in H_1$, $\rho(x,0) < \delta$. Then in particular
$$\sum_{i=1}^{a_N} \frac{|\langle x,h_i \rangle |}{2^i} < \delta ,$$
which implies $\|P_{F_N} x\| \leq \delta \cdot 2^{a_N}$. Now, for all $T\in\mathcal{T}$,
\begin{align*}
\rho(Tx,0)
&\leq \sum_{i=1}^{a_M}\frac{|\langle Tx,h_i \rangle |}{2^i} + \sum_{i>a_M}\frac{|\langle Tx,h_i \rangle |}{2^i} \\
&< \sum_{i=1}^{a_M}\frac{\|P_{F_M} Tx\|}{2^i} + \epsilon/4 \\
&< \|P_{F_M} Tx\|  + \epsilon/4 \\
&\leq \|P_{F_M} T P_{F_N} x \| + \|P_{F_M} T (I-P_{F_N}) x \| + \epsilon/4 \\
&< \delta \cdot 2^{a_N} + \|P_{F_M} T (I-P_{F_N}) x \| + \epsilon/4 \\
&\leq \|P_{F_M} T P_G(I-P_{F_N}) x \| + \|P_{F_M} T (I-P_G)(I-P_{F_N}) x \| + \epsilon/2 \\
&< \epsilon/4 + c + \epsilon/2 < \epsilon .
\end{align*}
All inequlities above follow from our choice of constants, and from noting that $(I-P_G)(I-P_{F_N}) x \perp G$.

(3) $\Rightarrow$ (1). We argue contrapositively. Let $V$ be a finite dimensional subspace of $H$, with an orthonormal basis $\{e_1,\ldots,e_n\}$. Assume that there is some $c>0$ and and infinite orthonormal sequence $\{e_k\}_{k>n}$ orthogonal to $\{e_1,\ldots,e_n\}$, such that
$$ \forall k>n . \exists T_k\in\mathcal{T} . \|P_V T_k e_k \| \geq c .$$
Without loss of generality, we may assume that $H = \mathrm{span}\{e_1, e_2, \ldots\}$ (set the elements of $\mathcal{T}$ to be zero on the complement and then throw it away). Define $\rho$ as we did in the proof of (2) $\Rightarrow$ (3). Then $\rho(e_k,0) \rightarrow 0$, but
\begin{align*}
\rho(T_k e_k,0)
&\geq \sum_{i=1}^{a_n} \frac{|\langle T_k e_k, h_i \rangle|}{2^i} \\
&\geq \sum_{i=1}^{n} \frac{|\langle T_k e_k, e_i \rangle|}{2^{a_n}} \\
&\geq \frac{1}{2^{a_n}}\sqrt{\sum_{i=1}^{n} |\langle T_k e_k, e_i \rangle|^2} \\
&\geq \frac{1}{2^{a_n}} \|P_V T_k e_k \| \geq \frac{c}{2^{a_n}} ,
\end{align*}
thus, $\mathcal{T}$ is not UEC at 0.
\end{proof}

\begin{corollary}
Let $\mathcal{U} \subset B_1$ be a family of isometries (or of maps uniformly bounded below). For every finite dimensional subspace $V \subset H$, put
$$\mathcal{U}^{-1}(V) = \{h\in H: \exists U \in \mathcal{U} . Uh \in V\} .$$
A necessary condition for $\mathcal{U}$ to be UEC on $H_1$ is
that for every finite dimensional subspace $V \subset H$,
    $$\mathrm{dim}\Big(\mathcal{U}^{-1}(V)\cap V^\perp\Big) < \infty .$$
\end{corollary}

\begin{corollary}\label{cor:suff}
Let $\mathcal{T} \subset B_1$ be a family of maps on $H_1$. For every finite dimensional subspace $V \subset H$, put
$$\mathcal{T}(V) = \{T h : h\in H , T \in \mathcal{T} \} .$$
A sufficient condition for $\mathcal{T}$ to be UEC on $H_1$ is that there exists an orthonormal basis $\{e_n\}_{n\in\mathbb{N}}$ and a number $K\in \mathbb{N}$ such that for all $N\in\mathbb{N}$,
\be\label{eq:assumption}
\mathcal{T}\big(\mathrm{span}\{e_{N+K}, e_{N+K+1}, \ldots \} \big) \perp \mathrm{span}\{e_1, \ldots, e_N\} .
\ee
\end{corollary}
\begin{proof}
For $K=0$, this is an immediate consequence of Theorem \ref{thm:eq}, (2) $\Rightarrow$ (3). Now let $K\geq 1$. Denote by $S_{r}$ the right shift with respect to the basis $\{e_n\}$. Then $\cT' : = S_r^K \cT  = \{S_r^K T \}_{T\in\cT}$ satisfies (\ref{eq:assumption}) with $K=0$, so $\cT = (S_r^* )^K \cT'$ must be UEC.
\end{proof}
Let us interpret Corollary \ref{cor:suff}. The corollary tells us that if there is an orthonormal basis $\cE = \{e_n\}_{n\in\mathbb{N}}$ and a $K\in \mathbb{N}$  such that for every $T \in \mathcal{T}$, the matrix of $T$ with respect to $\cE$ satisfies $t_{i,j} = 0$ for all elements above the $K$'th diagonal, then $\mathcal{T}$ is UEC. This seems like a condition one could check.

The statement of Proposition \ref{prop:Tstar} makes one wonder wether it possible to have $\mathcal{T}\subseteq B_1$ UEC
with $\mathcal{T}^*$ not UEC. The above condition makes it possible to draw an example.
\begin{expl}
Let $S_r$ be the right shift in $\ell^2(\mathbb{N})$. Since the family $\mathcal{T} = \{S_r^n\}_{n\in\mathbb{N}}$ clearly satisfies the above interpretation of corollary \ref{cor:suff}, this family is UEC. With $\mathcal{T}^*$ the situation is reversed (it maps infinitely many standard basis vectors isometrically into every finite dimensional space), thus it is not UEC.
\end{expl}
This example also shows that the assumption that $\mathcal{T}^*$ be UEC in Proposition \ref{prop:Tstar} cannot
be replaced by the assumption that $\mathcal{T}$ be UEC.

\bibliographystyle{amsplain}

\end{document}